\documentclass[12pt,letterpaper,reqno]{amsart}

\usepackage{times} 
\usepackage[T1]{fontenc} 
\usepackage{mathrsfs} 
\usepackage{latexsym} 
\usepackage[dvips]{graphics}
\usepackage{epsfig}
\usepackage{amsmath,amsfonts,amsthm,amssymb,amscd} 
\input amssym.def 
\input amssym.tex

\newcommand\be{\begin{equation}} 
\newcommand\ee{\end{equation}}
\newcommand\bea{\begin{eqnarray}} 
\newcommand\eea{\end{eqnarray}} 
\newcommand\bi{\begin{itemize}}
\newcommand\ei{\end{itemize}} 
\newcommand\ben{\begin{enumerate}} 
\newcommand\een{\end{enumerate}}
\newcommand\bc{\begin{center}} 
\newcommand\ec{\end{center}} 
\newcommand\ba{\begin{array}} 
\newcommand\ea{\end{array}}

\theoremstyle{definition} 

\begin{document}

\title[Limit points in the range of the commuting probability function on finite groups] 
{Limit points in the range of the commuting probability function on finite groups}


\author{Peter Hegarty} \address{Department of Mathematical Sciences, 
Chalmers University Of Technology and University of Gothenburg,
41296 Gothenburg, Sweden} \email{hegarty@chalmers.se}



\date{\today}

\begin{abstract}
If $G$ is a finite group, then Pr$(G)$ denotes the fraction of ordered pairs of elements of $G$ which commute. We show that if $l \in \left( \frac{2}{9}, 1 \right]$ is a limit point of the function Pr on finite groups, then $l \in \mathbb{Q}$ and there exists an $\epsilon = \epsilon_l > 0$ such that Pr$(G) \not\in \left( l - \epsilon_l, l \right)$ for any finite group $G$. 
These results lend support to some old conjectures of Keith Joseph.  
\end{abstract}


\maketitle

\setcounter{equation}{0}

\setcounter{equation}{0}

\setcounter{section}{0}


\section{Introduction}

Let $\mathcal{G}$ denote the family of all finite groups and define the function Pr$: \mathcal{G} \rightarrow \mathbb{Q} \cap (0,1]$ as follows{\footnote{Various alternative notations for this function appear in the literature, for example: cp$(G)$, mc$(G)$, $d(G)$.}}: for $G \in \mathcal{G}$, 
\be
{\hbox{Pr}}(G) := \frac{ \# \{(x,y) \in G \times G : xy = yx \}}{|G|^2}.
\ee
We call Pr$(G)$ the {\em commuting probability}{\footnote{Some authors, including for example Lescot, who has written a number of papers about this function (see \cite{L}), calls it the {\em commutativity degree} rather than the commuting probability.}} for $G$, in other words it is the probability that a uniformly random ordered pair of elements of $G$ commute. It is easy to see that Pr$(G) = k(G)/|G|$, where $k(G)$ denotes the number of conjugacy classes in $G$. There is quite a lot of literature on the properties of the function Pr. Much of this dates from a decade-or-so long period from the late 1960s to the late 1970s. In 1968, Erd\H{o}s and Tur\'{a}n \cite{ET} and, independently, Newman \cite{New}, proved that, for any finite group $G$, 
\be
{\hbox{Pr}}(G) \geq \frac{\log_2 \log_2 |G|}{|G|}.
\ee
This is a quantitative version of the classical fact that there are only finitely many finite groups with a given number of conjugacy classes. There have been various improvements on (1.2) since, as well as improved lower bounds for particular classes of groups: see \cite{Ke} and the references therein. It is easy to prove, as suggested by (1.2), that Pr$(G)$ can be arbitrarily close to zero. On the other hand, it is obvious that Pr$(G) = 1$ if and only if $G$ is abelian. Gustafson \cite{G} seems to have been the first to record, in 1973, the fact that, if $G$ is non-abelian, then Pr$(G) \leq 5/8$, with equality being achieved if and only if $G/Z(G) \cong C_2 \times C_2$. The intervening years have witnessed a minor flood of papers whose common theme is to show that, if Pr$(G)$ is sufficiently large, then $G$ is ``close to abelian'' in some well-defined structural sense. There are basically two types of results here:
\par (i) Rusin \cite{R} completely classifies all finite groups for which Pr$(G) > 11/32$. A recent work of Das and Nath \cite{DN} does the same{\footnote{There were some errors in Rusin's paper, which were pointed out and corrected by Das-Nath. See also the proof of Lemma 2.1 below.}} for all odd order groups satisfying Pr$(G) \geq 11/75$. The point here is that one can completely determine all finite groups $G$ for which Pr$(G)$ lies above a certain threshold.     
\par (ii) There are also results of a more general character, whose purpose is to determine some general ``abelian-like feature'' of a group $G$ for which Pr$(G)$ is bounded away from zero. A fairly recent work of Guralnick and Robinson \cite{GR} contains a number of impressive such results. In particular, they prove (\cite{GR}, Theorem 10) that 
\be
{\hbox{Pr}}(G) \leq (G:F)^{-1/2},
\ee
where $F$ is the Fitting subgroup of $G$. Hence, if Pr$(G)$ is bounded away from zero, then $G$ contains a normal, nilpotent subgroup of bounded index. A similar, but even more striking result was shown much earlier by Peter Neumann:
\\
\\
{\bf Theorem 1.1 (P.M. Neumann \cite{Neu}).} {\em For each $\epsilon > 0$, there exist positive integers $n_{1}(\epsilon), n_{2}(\epsilon)$ such that the following holds: If $G$ is a finite group satisfying Pr$(G) > \epsilon$, then $G$ possesses normal subgroups $H,K$ such that
\par (i) $K \subseteq H$,
\par (ii) $|K| \leq n_{1}(\epsilon)$,
\par (iii) $(G:H) \leq n_{2}(\epsilon)$,
\par (iv) $H/K$ is abelian.} 
\\
\\
In this paper we are basically interested in the following question:
\\
\par {\em What is the range of the function Pr inside $\mathbb{Q} \cap (0,1]$ ?}
\\
\\
Gustafson's observation makes this question very natural, as it implies that the range is not the whole of $\mathbb{Q} \cap (0,1]$. This is in stunning contrast to the situation for semigroups: see \cite{PS}. In the case of groups, further insight into the range of Pr can be gleaned, if indirectly, from the various papers cited above. As far as we know, however, only one person ever seems to have seriously considered general structural properties of Pr$(\mathcal{G})$ as a set. Keith Joseph wrote his Ph.D. thesis \cite{J1} in 1969 on the commuting probability, but it was never published. Eight years later, however, in a short note in the Monthly \cite{J2}, he posed three very interesting conjectures about the range of Pr, which we reproduce below{\footnote{The second and third conjectures were stated slightly differently by Joseph.}}. Note that, by a {\em limit point} of a set $S \subseteq \mathbb{R}$, we mean an element $l \in \mathbb{R}$ such that there is a sequence $(s_n)$ of elements of $S$ tending to $l$. In particular, every element of $S$ itself is a limit point of $S$. 
\\
\\
{\bf Joseph's first conjecture.} {\em Every limit point of Pr$(\mathcal{G})$ is rational.}
\\
{\bf Joseph's second conjecture.} {\em If $l$ is a limit point of Pr$(\mathcal{G})$, then there exists $\epsilon = \epsilon_l > 0$ such that Pr$(\mathcal{G}) \cap (l-\epsilon_l,l) = \phi$.}
\\
{\bf Joseph's third conjecture.} {\em Pr$(\mathcal{G}) \cup \{0\}$ is a closed subset of $\mathbb{R}$.}
\\
\\
Unfortunately, and despite the large amount of literature in existence today on the commuting probability, Joseph's conjectures seem to have suffered the same fate as the rest of his work and been thoroughly neglected{\footnote{As already mentioned, his Ph.D. thesis was never published. Neither is he listed as an author on MathSciNet. His Monthly article appeared as a ``research problem'', not as a regular paper.}}. There seems to have been essentially no progress on any of them, and we found only a handful of references to his work. For example, Rusin \cite{R} refers to Joseph's questions at the end of his paper. He notes that, since his methods rely heavily on the elementary estimate
\be
{\hbox{Pr}}(G) \leq \frac{1}{4} + \frac{3}{4} \frac{1}{|G^{\prime}|},
\ee
they cannot, even in principle, yield any insight into the range of Pr inside $(0,1/4]$. Rusin wonders explicitly whether the range is dense inside $(1/5,1/4)$. In a recent Master's thesis, Castelaz{\footnote{who has since married and appears on MathSciNet as Anna Keaton, see: \par https://mthsc.clemson.edu/directory/view$\underline{}$person.py?person$\underline{}$id=198}} (\cite{C}, end of Chapter 5) essentially repeats Joseph's questions, but does not provide any new insights. Only the following facts seem to be known, and all but the last appear to have been already known to Joseph. For proofs of (i)-(iv), see \cite{C}:
\\
\\
(i) zero is a limit point, not in the range of Pr.
\\
(ii) for each prime $p$, $1/p$ is a limit point of Pr$(\mathcal{G})$ and there exists a group $G$ with Pr$(G) = 1/p$.
\\
(iii) since the function Pr is multiplicative, i.e.:
\be
{\hbox{Pr}}(G_1 \times G_2) = {\hbox{Pr}}(G_1) \times {\hbox{Pr}}(G_2),
\ee
it follows that the same is true as in (ii) for every number $a/n \in (0,1]$, such that
$n \in \mathbb{N}$ and $a \in {\hbox{Pr}}(\mathcal{G})$. 
\\
(iv) Joseph's second conjecture is known to hold for $q = 1/2$ (with $\epsilon_{1/2} = 1/16$). 
\\
(v) The following is a corollary of a result of Gavioli et al:
\\
\\
{\bf Theorem 1.2 (Corollary to Theorem 3 in \cite{GMMPS}).} {\em Let $p \geq 7$ be a prime and let $G$ be a $p$-group. Then 
\be
{\hbox{Pr}}(G) \not\in \left[ \frac{5}{p^4}, \frac{1}{p^3} \right].
\ee}
In this paper we shall provide some new evidence for Joseph's first two conjectures by proving the following result:
\\
\\
{\bf Main Theorem.} {\em If $l \in \left( \frac{2}{9}, 1 \right]$ is a limit point of the set Pr$(\mathcal{G})$, then 
\par (i) $l \in \mathbb{Q}$,
\par (ii) there exists an $\epsilon = \epsilon_l > 0$ such that Pr$(\mathcal{G}) \cap \left( l - \epsilon_l, l \right) = \phi$.} 
\\
\\
Note that, as stated in (iv) above, the second assertion in our theorem is already known for $l = 1/2$. As well as giving a bit of new evidence in support of Joseph's long-dormant program, our result provides the first insight into the range of Pr below $1/4$, for arbitrary groups (Theorem 1.2 gives some insight for $p$-groups and the results in \cite{DN} for odd order groups), as sought by Rusin. After some preliminaries in Section 2, we will prove our theorem in Section 3. In Section 4 we will discuss the prospects for extending our methods in the hope of answering Joseph's questions definitively. 

\setcounter{equation}{0}

\section{Preliminaries}

Throughout the rest of this paper, all groups are finite unless explicitly stated otherwise. Recall that groups $G_1$ and $G_2$ are said to be {\em isoclinic} if there are isomorphisms 
\be
\phi: G_1 / Z_1 \rightarrow G_2 /Z_2, \;\;\; \psi: G_{1}^{\prime} \rightarrow G_{2}^{\prime} \;\;\;\;\; (Z_1 = Z(G_1), \; Z_2 = Z(G_2)),
\ee
such that, for all $x,y \in G_1$, 
\be
\psi \left( [Z_1 x, Z_1 y] \right) = \left[ \phi (Z_1 x), \phi (Z_1 y) \right].
\ee
Isoclinism is an equivalence relation on the set of all finite groups and each equivalence class contains a so-called {\em stem group}, i.e.: a group $G$ such that $Z(G) \subseteq G^{\prime}$.
\\
\\
First, we collect a number of facts about the function Pr:
\\
\\
{\bf Lemma 2.1.} {\em (i) If $G$ is a non-abelian group then Pr$(G) \leq 5/8$. Moreover, if Pr$(G) > 11/32$ then one of the following must hold :
\par (a) $|G^{\prime}| = 2$, $G/Z(G) \cong C_{2}^{2s}$ for some $s \in \mathbb{N}$ and Pr$(G) = \frac{1}{2} \cdot (1 + 2^{-2s})$.
\par (b) $G^{\prime} \cong C_3$, $G/Z(G) \cong S_3$ and Pr$(G) = 1/2$.
\par (c) Pr$(G) \leq 7/16$ and $|G/Z(G)| \leq 16$.
\\
(ii) Let $p$ be a prime and $G$ a $p$-group satisfying $G^{\prime} \subseteq 
Z(G)$. Then
\be
{\hbox{Pr}}(G) = \frac{1}{|G^{\prime}|} \left( 1 + \sum_{K} \frac{(p-1) (G^{\prime}:K)}{p^{s(K) - 1}} \right),
\ee
where the sum is taken over all subgroups $K$ of $G^{\prime}$ for which $G^{\prime}/K$ is cyclic and the integer $n(K)$ is defined by 
\be
p^{s(K)} \stackrel{{\hbox{def}}}{=} \frac{|G|}{| \{x \in G : [G,x] \subseteq K\}|}.
\ee
In particular, if $G^{\prime} \cong C_p$ and $G/Z(G) \cong C_{p}^{2s}$, then
\be
{\hbox{Pr}}(G) = \frac{1}{p} \left( 1 + \frac{p-1}{p^{2s}} \right).
\ee
(iii) If $|G^{\prime}| = 4$ and $|G^{\prime} \cap Z(G)| = 2$ then 
\be
{\hbox{Pr}}(G) = \frac{1}{4} \left( 1 + \frac{1}{4} + \frac{1}{2^{2s+1}} \right),
\ee
where $2^{2s} = [C_{G}(G^{\prime}):Z(C_{G}(G^{\prime}))]$.
\\
(iv) If $G^{\prime} \cong C_6$ and $G^{\prime} \cap Z(G) \cong C_2$ then 
Pr$(G) = 1/4 + 1/2^s$ for some $s \geq 3$.
\\
(v) For any fixed group $K$, there are only finitely many groups $G$ satisfying $G^{\prime} \cong K$ and $Z(G) = \{1\}$. In particular, there are only finitely many possibilities for Pr$(G)$ under these conditions. 
\\
(vi) Let $G$ be a non-abelian group and let $d$ be the minimum degree of a non-linear irreducible representation of $G$. Then 
\be
\frac{1}{|G^{\prime}|} < {\hbox{Pr}}(G) \leq \frac{1}{d^2} + \left( 1 - \frac{1}{d^2} \right) \frac{1}{|G^{\prime}|}.
\ee
(vii) If $H$ is a subgroup of $G$ then Pr$(H) \geq {\hbox{Pr}}(G)$.
\\
(viii) If $N$ is a normal subgroup of $G$, then Pr$(G) \leq {\hbox{Pr}}(N) \cdot {\hbox{Pr}}(G/N)$.
\\
(ix) Let $N$ be a normal subgroup of $G$. Let $c$ be an integer such that
every subgroup of $G/N$ contains at most $c$ conjugacy classes. Let $k_{G}(N)$ 
denote the number of $G$-orbits of $N$, when $G$ acts on $N$ by conjugation. Then $k(G) \leq c \cdot k_{G}(N)$ and hence 
\be
{\hbox{Pr}}(G) \leq \frac{c}{|G/N|} \cdot \frac{k_{G}(N)}{|N|}.
\ee
(x) If the groups $G_1$ and $G_2$ are isoclinic, then Pr$(G_1) = {\hbox{Pr}}(G_2)$. 
\\
(xi) Let $n \in \mathbb{N}$ and let $\mathcal{Z}_n$ denote the collection of all finite groups $G$ for which \\ $(G:Z(G)) \leq n$. Then Pr$(\mathcal{Z}_n)$ is a finite set.}
\\
\\
{\sc Proof}: Parts (i),(ii),(iv) and (v) can be found in \cite{R}. Part (iii) is proven in \cite{DN}, and they correct an erroneous form of (2.6) which appeared in \cite{R}. Parts (vi)-(ix) can be found in \cite{GR} and part (x) in \cite{L}. We have not seen part (xi) written down anywhere (though similar statements appear in \cite{NiS}), so we give the easy proof: Let $n$ be given and pick $G$ such that $(G:Z(G)) \leq n$. Then $|G'|$ is bounded in terms of $n$, by Schur's classical result \cite{S}. Now, by part (vi), we may assume $G$ is a stem-group, in which case it follows that $|Z(G)|$ is also bounded. Thus $|G|$ is bounded, and hence there are only finitely many possible values of Pr$(G)$.
\\
\\
{\bf Lemma 2.2.} {\em For elements $x,y,z,w$ in a group $G$ one has 
\be
[xy,zw] = [x,w]^{y} [x,z]^{wy} [y,w] [y,z]^{w}.
\ee}
{\sc Proof:} Simple exercise.
\\
\\
{\bf Lemma 2.3.} {\em Let $G$ be a finite subgroup of {\hbox{PGL}}$(2,\mathbb{C})$. Then $G$ is isomorphic to one of the following:
\par (i) a cyclic group $C_n$, for some $n \in \mathbb{N}$,
\par (ii) a dihedral group $D_n$ of order $2n$, for some $n \geq 2$,
\par (iii) the tetrahedral group $A_4$,
\par (iv) the octahedral group $S_4$,
\par (v) the icosahedral group $A_5$.
\\
\\
Moroever, we have that 
\be
{\hbox{Pr}}(D_n) = \left\{ \begin{array}{lr} \frac{n+6}{4n}, & {\hbox{if $n$ is even}}, \\ \frac{n+3}{4n}, & {\hbox{if $n$ is odd,}} \end{array} \right. \;\;\; {\hbox{Pr}}(A_4) = \frac{1}{3}, \;\;\; {\hbox{Pr}}(S_4) = \frac{5}{24}, \;\;\; {\hbox{Pr}}(A_5) = \frac{1}{12}.
\ee}
{\sc Proof:} The classification of the finite subgroups of PGL$(2,\mathbb{C})$ is classical; see, for example, \cite{Kl}. The values of Pr for these groups can be easily verified. 
\\
\\
Before stating our last preliminary result, let us introduce some terminology which will simplify the statement of this and succeeding results:
\\
\\
{\bf Definition 2.1.} A subset $S \subseteq \mathbb{R}$ will be called {\em good} if, for every limit point $l$ of $S$, the following hold:
\par (i) $l \in \mathbb{Q}$,
\par (ii) there exists $\epsilon = \epsilon_l > 0$ such that $S \cap (l-\epsilon_l,l) = \phi$. 
\\
\\
Observe that any subset of a union of finitely many good sets is also good. 
\\
\\
{\bf Lemma 2.4.} {\em Let $n$ be a fixed positive integer and let 
\be
\mathcal{S}_{n} := \left\{ \sum_{i=1}^{n} \frac{1}{x_i} : x_i \in \mathbb{N} \right\}.
\ee
Then Pr$(\mathcal{S}_n)$ is a good set.}
\\
\\
{\sc Proof:} It is a classical fact that, for every fixed $n \in \mathbb{N}$ and $q \in \mathbb{Q}_{+}$, the equation 
\be
\sum_{i=1}^{n} \frac{1}{x_i} = q
\ee
has only finitely many positive integer solutions $(x_1,...,x_n)$. If one examines the standard proof of this fact (which is basically just an induction on 
$n$), one easily sees that it in fact implies what is claimed in the lemma.

\setcounter{equation}{0}

\section{Proof of Main Theorem}

The following lemma is the crucial ingredient in our proof, and we have not seen it written down before. The reader should observe the connection to Theorem 1.1, more about which will be said in the next section. 
\\
\\
{\bf Lemma 3.1.} {\em Let $n \in \mathbb{N}$ and let $\mathcal{A}_n$ be the collection of all finite groups possessing a normal, abelian subgroup of index $n$.
Then Pr$(\mathcal{A}_n)$ is a good set.}
\\
\\
{\sc Proof:} Fix $n \in \mathbb{N}$, a group $G$ and a normal abelian subgroup $H$ such that $(G:H) = n$. Let 
\be
G = \bigsqcup_{i=1}^{n} Hx_i, \;\;\;\; (x_1 = 1)
\ee
be a decomposition of $G$ into cosets of $H$. For each ordered pair $(i,j)$ of indices from the set $\{1,...,n\}$, let 
\be
S_{ij} := \{ (h_1,h_2) \in H \times H : [h_1 x_i, h_2 x_j] = 1 \}.
\ee
Thus
\be
{\hbox{Pr}}(G) = \frac{1}{|G|^2} \cdot \sum_{i,j = 1}^{n} |S_{ij}|.
\ee
Since $H$ is abelian, it follows easily from Lemma 2.2 that, for any $g \in G$ the map $h \mapsto [h,g]$ is an endomorphism of $H$, whose kernel is $C_{H}(g)$. Let $H_g := [H,g]$ and $n_g := (H:C_{H}(g))$. Thus $H_g$ is a subgroup of $H$ of order $n_g$. For each $i = 1,...,n$ above, let $H_i := H_{x_i}$ and $n_i := n_{x_i}$. For each ordered pair $(i,j)$, set $H_{ij} := H_i \cap H_j$ and $n_{ij} := |H_{ij}|$. Thus $n_{ij}$ is a common divisor of $n_i$ and $n_j$. Now fix a pair $(i,j)$ and set $h_{ij} := [x_j,x_i]$. This is a fixed element of $H$. If $h_1,h_2 \in H$ then, using Lemma 2.2, it is easy to check that 
\be
[h_1 x_i, h_2 x_j] = 1 \Leftrightarrow [h_{1}^{x},y] = h_{ij} [h_{2}^{y},x].
\ee
Let 
\be
\hat{H}_{ij} := H_j \cap (h_{ij} H_i) = \{ h \in H_j : h = h_{ij}u, \;\; {\hbox{for some $u \in H_i$}} \}.
\ee
It's easy to see that either $\hat{H}_{ij} = \phi$ or is a single coset in $H$ of
the subgroup $H_{ij}$. In the former case, the right-hand side of (3.4) has no solutions. In the latter case, we can count the number of solutions as follows: first, we pick $h_1$ such that $[h_{1}^{x},y] \in \hat{H}_{ij}$. The number of possible choices is just $|H| \cdot \frac{|H_{ij}|}{|H_j|} = |H| \cdot \frac{n_{ij}}{n_j}$. Having chosen $h_1$, we pick $h_2$ so that the right-hand equation in (3.4) is satisfied. The number of choices for $h_2$ is just $(H:H_i) = |H| \cdot \frac{1}{n_i}$. Summarising, we have shown that
\be
|S_{ij}| = \left\{ \begin{array}{lr} 0, & {\hbox{if $\hat{H}_{ij} = \phi$}}, \\
|H|^{2} \cdot \frac{n_{ij}}{n_i n_j}, & {\hbox{otherwise}}. \end{array} \right.
\ee
Hence, the expression (3.3) for Pr$(G)$ has the form 
\be
{\hbox{Pr}}(G) = \frac{1}{n^2} \cdot \sum_{k=1}^{L} \frac{1}{x_k},
\ee
where each $x_k$ is a positive integer, $x_1 = 1$ and $1 \leq L \leq n^2$. Here, $x_1$ corresponds to the term $S_{11}$ and the fact that $L$ may be less than $n^2$ corresponds to the fact that some of the $S_{ij}$ may be empty. Further, note that the numbers $x_k$ are not independent of one another, since the same can be said of the numbers $|S_{ij}|$. However, this just makes our life easier. We conclude that the set of possible values for Pr$(G)$, in the notation of Lemma 2.4, is contained in the set 
\be
\frac{1}{n^2} \cdot \left( \bigcup_{k=1}^{n^2} \mathcal{S}_k \right).
\ee
Lemma 2.4 thus directly implies the claim of Lemma 3.1.
\\
\\
{\bf Corollary 3.2.} {\em For each $n \in \mathbb{N}$, let $\mathcal{A}^n$ denote the collection of all finite groups possessing an abelian subgroup of index at most $n$. Then Pr$(\mathcal{A}^n)$ is a good set.}
\\
\\
{\sc Proof:} If $A$ is an abelian subgroup of $G$ of index at most $n$, then Core$_{G}(A)$ is an abelian, normal subgroup of index at most $n!$. Hence $\mathcal{A}^{n} \subseteq \cup_{k=1}^{n!} \mathcal{A}_k$ and we can apply Lemma 3.1.
\\
\\
We are now ready to prove the Main Theorem in a sequence of steps. Let $G$ be a
non-abelian group satisfying Pr$(G) > 2/9$. By Lemma 2.1(x), we may assume that $Z(G) \subseteq G^{\prime}$. 
\\
\\
{\em Step 1:} Let $d$ be the minimum degree of a non-linear irreducible representation of $G$. From Lemma 2.1(vi) we deduce that either $|G^{\prime}| < 8$ or 
$d = 2$.
\\
\\
{\em Step 2:} First suppose $|G'| < 8$. Since $Z(G) \subseteq G^{\prime}$, it is easy to check that either $G$ is covered by parts (iii), (iv) and (v) of Lemma 2.1, or $G$ is nilpotent and a direct product of $p$-groups satisfying (2.5). It's then just a matter of verifying that the Main Theorem is satisfied in these cases.    
\\
\\
{\em Step 3:} So we may suppose $d=2$. Let $\phi$ be an irreducible representation of $G$ of degree $2$. Let $\pi: {\hbox{GL}}(2,\mathbb{C}) \rightarrow {\hbox{PGL}}(2,\mathbb{C})$ be the natural projection and set 
$K := {\hbox{ker}}(\phi)$, $L := {\hbox{ker}}(\pi \circ \phi)$. Then $G/L$ is isomorphic to a finite subgroup of PGL$(2,\mathbb{C})$, hence to one of the non-cyclic groups listed in Lemma 2.2. 
\par First suppose that $G/L \cong A_4, S_4$ or $A_5$. Since Pr$(G) > 2/9$, 
the second and third options are immediately ruled out by (2.10) and part (viii) of Lemma 2.1. In the case of $A_4$, the same analysis, together with Lemma 2.1(i), implies that $K$ must be abelian. But then we can apply Lemma 3.1, and the Main Theorem is satisfied. 
\par So we may suppose that 
\be
G/K \cong Z \cdot 2D_n, \; {\hbox{for some $n \geq 2$, where $Z$ is a finite cyclic group}}.
\ee
First suppose $n \geq 3$. Then $G/L \cong D_n$, say 
\be
G/L = \; <La, Lb | a^{n}, b^2, (ab)^2 \in L>.
\ee
We now consider two separate cases:
\\
\\
{\sc Case 1:} $n \geq 15$. 
\\
\\
Since $\frac{5}{8} \left( \frac{n+6}{4n} \right) < \frac{2}{9}$ for all $n \geq 15$, we can argue as before that $L$ must be abelian. By Lemma 2.1(ix),
\be
\frac{c}{2n} \frac{k_{G}(L)}{|L|} > \frac{2}{9},
\ee
where the number $c$ is such that every subgroup of $D_n$ contains at most $c$ conjugacy classes. Clearly we can take $c = n$, whence (3.11) becomes 
\be
\frac{k_{G}(L)}{|L|} > \frac{4}{9}.
\ee
In other words, the average size of a $G$-orbit in $L$ is less than $9/4$. 
It follows easily that $a^k \in C_{G}(L)$ for some $k = O(1)$, independent of $n$. Let $N := \; <L,a^k>$. Then $N$ is an abelian, normal subgroup of $G$ of bounded index, so the Main Theorem holds, by Lemma 3.1.
\\
\\
{\sc Case 2:} $3 \leq n \leq 14$.
\\
\\
Then $|G/L|$ is bounded. If $(L:Z(L))$ were also bounded, then $Z(L)$ would be a normal, abelian subgroup of $G$ of bounded index and we could apply Lemma 3.1 again. So we may suppose $L$ is non-abelian. Since Pr$(G/L) \leq \frac{5}{8}$ and since $\frac{11}{32} \times \frac{5}{8} < \frac{2}{9}$, Lemma 2.1(i) would still imply that $(L:Z(L))$ were bounded, unless $|L^{\prime}| = 2$ and $L/Z(L) \cong C_{2}^{2s}$ for some $s \in \mathbb{N}$. Thus, $L^{\prime} \subseteq Z(G)$. 
We can still apply Lemma 2.1(ix) to conclude that
\be
\frac{k_{G}(L)}{|L|} > \frac{8}{5} \times \frac{2}{9} = \frac{16}{45}.
\ee
In other words, the average size of a $G$-orbit in $L$ is less than $45/16$. This must imply that $(L:L_1)$ is bounded, where $L_1 = \{x \in L : [G,x] \subseteq L^{\prime} \}$. Now $G^{\prime}/L^{\prime} \cong (G/L_1)^{\prime}$. By Lemma 2.1(i), either $|G^{\prime}| \leq 6$ or Pr$(G/L_1) \leq \frac{7}{16}$. The first alternative takes us back to {\em Step 2}. From the second alternative and Lemma 2.1(viii) we conclude that Pr$(L_1) > \frac{32}{63} = \frac{1}{2} + \frac{1}{126}$ and hence, by Lemma 2.1(i), that $(L_1 : Z(L_1))$ is bounded. But since $(G:L)$ and $(L:L_1)$ are also bounded, we conclude that $(G:Z(L_1))$ is bounded and we can apply Lemma 3.1 one more time to conclude the analysis of {\sc Case 2}. 
\\
\\
We are now left with the possibility that $n = 2$ in (3.9). Set $G/K := Q$. Then $|Q^{\prime}| = 2$, $Q/Z(Q) \cong C_2 \times C_2$ and $Z(Q)$ is a cyclic group. A priori, the order of $Q$ may be unbounded, but the crucial thing is that $Q$ has a cyclic subgroup of bounded index. Let $q \in G$ be such that $Kq$ generates $Z(Q)$. Since Pr$(Q) = 5/8$ we can repeat the analysis from {\sc Case 2} above to conclude that either $|G^{\prime}| \leq 6$ or $K$ contains an abelian subgroup $K_2$ of bounded index such that $K_2 \lhd G$. So we may suppose the latter holds. Lemma 2.1(ix) still applies and, as in (3.13), we have that the average size of a $G$-orbit in $K_2$ is less than $45/16$. Let $K_3 := \{x \in K_2 : (G:C_{G}(x)) \leq 2\}$. A priori, $K_3$ may not be a subgroup of $K_2$, however we must have that $|K_3|/|K_2|$ is bounded away from zero. In addition, since $(G:C_{G}(x)) \leq 2$ for all $x \in K_3$, it follows that $K_3 \subseteq C_{g}(q^2)$. Let 
$K_4 := C_{G}(q^2) \cap K_2$ and $A := \; <K_4,q^2>$. Then $A$ is an abelian subgroup of $G$ of bounded index, and hence we can apply Corollary 3.2. 
This completes the proof of the Main Theorem. 
   
\setcounter{equation}{0}

\section{Discussion}

It is not true that if Pr$(G)$ is bounded away from zero, then $G$ contains a (normal) abelian subgroup of bounded index. Indeed, by Lemma 2.1(i), we see that this already fails for groups satisfying Pr$(G) > 1/2$. So we cannot prove Joseph's first two conjectures simply by using Lemma 3.1. Indeed, Theorem 1.1 seems to give the strongest possible structural result about groups for which Pr$(G)$ is bounded away from zero. Note that the structure described there includes the case when $|G^{\prime}|$ is bounded. Indeed, the strategy of our proof in the previous section began by appealing to Lemma 2.1(vi), which says that if Pr$(G)$ is bounded away from zero, then either $|G^{\prime}|$ or the minimum degree of a non-linear irreducible representation of $G$ is bounded. In the latter case, it is also interesting that a classical result of Jordan (see \cite{I}, Theorem 14.12) says that there is a function $f: \mathbb{N} \rightarrow \mathbb{N}$ such that a finite subgroup of GL$(n,\mathbb{C})$ must contain an abelian, normal subgroup of index at most $f(n)$. However, it is unlikely that we can say anything more than Theorem 1.1 about the structure of $G$ in general{\footnote{Neumann's proof does not use any representation theory, though he does make use of the well-known fact that if the sizes of the conjugacy classes in a group are bounded, then so is the size of the full commutator subgroup.}}. Therefore, it seems a crucial step in the analysis of Joseph's conjectures is to see if the sets Pr$(\mathcal{C}_n)$ are good, where $\mathcal{C}_n$ is the collection of all finite groups $G$ for which $|G^{\prime}| \leq n$. It is still not obvious to us how one would get from there and Lemma 3.1 to a full proof of Joseph's first two conjectures, but at least we have provided a possible roadmap.    
\par Finally, we have not said anything in this paper about Joseph's third conjecture, which seems more mysterious to us. 

\section*{Acknowledgement}

I thank Des MacHale for helpful discussions and for drawing my attention to 
several of the papers in the bibliography below.

\end{document}